\documentclass[10pt]{amsart}

\usepackage[T1]{fontenc}
\usepackage{lmodern}
\usepackage{microtype}
\usepackage[a4paper,margin=30mm]{geometry}
\usepackage{mathtools,amssymb}
\usepackage{enumitem}
\usepackage{needspace}
\usepackage[hidelinks]{hyperref}
\hypersetup{
  pdftitle={Negative normal degree and holomorphic sectional curvature pinching},
  pdfauthor={Xingluan Wang}
}

\newtheorem{theorem}{Theorem}[section]
\newtheorem{lemma}[theorem]{Lemma}
\newtheorem{corollary}[theorem]{Corollary}
\theoremstyle{remark}
\newtheorem{remark}[theorem]{Remark}
\numberwithin{equation}{section}

\newcommand{\CP}{\mathbb{CP}}
\newcommand{\OO}{\mathcal{O}}
\newcommand{\FF}{\mathbb{F}}

\DeclareMathOperator{\tr}{tr}

\title[Negative normal degree and holomorphic sectional curvature pinching]
      {Negative normal degree and holomorphic sectional curvature pinching}
\author{Xingluan Wang}
\thanks{E-mail: wxlmath@mail.ustc.edu.cn. School of Mathematical Sciences,
University of Science and Technology of China, Hefei 230026, P. R. China.}
\date{}

\begin{document}

\begin{abstract}
We obtain an upper bound for the local holomorphic sectional curvature pinching of a K\"ahler manifold containing a smooth rational curve of negative total normal degree. As an application, we prove that a compact K\"ahler surface with positive holomorphic sectional curvature and local pinching constant at least $1/9$ is biholomorphic to $\mathbb{CP}^2$, $\mathbb{CP}^1\times\mathbb{CP}^1$, or the first Hirzebruch surface $\mathbb F_1$. For pinching strictly greater than $1/9$, only the first two possibilities occur, confirming a conjecture of Yang and Zheng. We also determine the optimal local pinching constants on $\mathbb F_1$ and $\mathbb F_2$.
\end{abstract}

\maketitle

\section{Introduction and main results}\label{introduction}

\subsection{Background and motivation}
Let $(M^n,g)$ be a K\"ahler manifold of complex dimension $n\ge 2$ with positive holomorphic sectional curvature $H_g(v)=R(v,\bar v,v,\bar v)/|v|_g^4$ for $v\in T^{1,0}M\setminus\{0\}$.
For $p\in M$, write
\[
H_{\min}(p)=\min_{|v|_g=1}H_g(v),\qquad
H_{\max}(p)=\max_{|v|_g=1}H_g(v).
\]
For any nonempty subset $A\subset M$, define its local holomorphic pinching
constant by
\[
\lambda_A(g):=\inf_{p\in A}\frac{H_{\min}(p)}{H_{\max}(p)}.
\]
In particular, $\lambda_M(g)$ is the local pinching constant of $(M,g)$.
For comparison, its global pinching constant is
$\lambda_M^{\mathrm{glob}}(g):=\inf_M H_{\min}/\sup_M H_{\max}$. If $(M,g)$ is complete,
positive global pinching implies compactness by Tsukamoto's theorem
\cite{Tsukamoto}.

Classical global pinching results provide rigidity at relatively high thresholds. Kobayashi \cite{Kobayashi} obtained the homotopy type of $\CP^n$ under pinching greater than $11/13$. Bishop and Goldberg \cite{BG} improved the threshold to $4/5$. In complex dimension two, Greene and Pittie \cite{GP} showed that global pinching greater than $1/2$ forces the complex structure of $\CP^2$. Cao and Yang \cite{CY} subsequently established a gap theorem for compact K\"ahler manifolds with global pinching close to $1/2$. For local pinching, Yang and Zheng \cite[Proposition~2.6]{YZ} proved that a compact K\"ahler manifold with $\lambda_M(g)>1/2$ is biholomorphic to $\CP^n$, and classified the equality case. The relevance of $1/2$ follows from Berger's inequality  \cite{Berger}. For any \(p\in M\) and any orthogonal unit vectors \(u,v\in T_p^{1,0}M\),
\begin{equation}\label{b1}
R(u,\bar u,v,\bar v)
\ge
\left(\lambda_M(g)-\frac12\right)
H_{\max}(p).
\end{equation}
Thus, $\lambda_M(g)> 1/2$ implies positive orthogonal bisectional curvature. The Kähler\textendash{}Ricci flow produces a metric with positive bisectional curvature, and thus \(M\) is biholomorphic to \(\mathbb{CP}^n\). When \(\lambda_M(g)=1/2\), the classification follows from Gu\textendash{}Zhang's splitting theorem \cite[Theorem~1.3]{GuZhang}.

Below $1/2$, this nonnegativity need not hold. The first Hirzebruch surface \(\mathbb F_1\) admits a Hitchin-type Kähler metric with local pinching constant \(1/9\), while the flag threefold admits a canonical Kähler\textendash{}Einstein metric with local pinching constant \(1/4\). Yang and Zheng therefore asked whether a compact K\"ahler surface with $\lambda_M(g)>1/9$ must be biholomorphic to $\CP^2$ or $\CP^1\times\CP^1$ \cite[Conjecture~1.7]{YZ}. They also conjectured that if a compact Kähler threefold with $H > 0$ has its local holomorphic pinching constant $\lambda>1/4$, then it must be biholomorphic to a compact Hermitian symmetric space \cite[Conjecture~1.11]{YZ}.

\subsection{Main results}

Let $C\simeq\CP^1$ be a smoothly embedded holomorphic rational curve in $(M^n,g)$. Its holomorphic normal bundle is $N_{C/M}:=T^{1,0}M|_C/T^{1,0}C$.
In this paper, we derive an upper bound for \(\lambda_M(g)\) in terms of the negative normal degree \(d\) and the integral \(S\) of the squared norm of the second fundamental form \(\Pi\) of \(C\). Here $d$ and $S$ are defined by 
\[
d:=-\deg N_{C/M}, \qquad S:=\int_C|\Pi(e,e)|_g^2\,\omega_C,
\]
where $e$ is any local unit tangent
vector along $C$, and $\omega_C=\omega_g|_C$. For $d>0$, define
\[
\mu(n,d,S):=\frac{2\pi d+S}{(n-1)(4\pi+S)},\qquad
\Phi(\mu):=
\begin{cases}
\frac12-\mu,&0<\mu\le\frac14,\\[2pt]
(1+4\mu)^{-2},&\mu\ge\frac14.
\end{cases}
\]
The following estimate applies in every complex dimension $n\ge2$.

\begin{theorem}\label{main1}
Let $(M^n,g)$ be a K\"ahler manifold with positive holomorphic sectional curvature. If $M$ contains a smooth rational curve $C$ with $\deg N_{C/M}=-d<0$, then
\begin{equation}\label{main-bound}
\lambda_M(g)\le\lambda_C(g)
\le\Phi\bigl(\mu(n,d,S)\bigr)
\le\Phi\left(\frac{\min\{d/2,1\}}{n-1}\right).
\end{equation}
\end{theorem}

In dimension two,
a $(-1)$-curve gives the upper bound $1/9$. This observation leads to the following classification.

\Needspace{11\baselineskip}
\begin{theorem}\label{main2}
Let $(M^2,g)$ be a compact connected K\"ahler surface with positive holomorphic sectional curvature. Then:
\begin{enumerate}[label=\textup{(\roman*)},leftmargin=*,itemsep=1pt]
\item If $\lambda_M(g)>1/9$, then $M$ is biholomorphic to $\CP^2$ or $\CP^1\times\CP^1$.
\item If $\lambda_M(g)\ge1/9$, then $M$ is biholomorphic to $\CP^2$, $\CP^1\times\CP^1$, or $\FF_1$.
\end{enumerate}
\end{theorem}

\begin{remark}
In complex dimension two, when $\lambda_M(g)>1/9$, Theorem~\ref{main1} excludes rational curves of negative normal degree. At the critical value $\lambda_M(g)=1/9$, such curves may still exist, but they must be totally geodesic $(-1)$-curves.
\end{remark}

The total normal degree does not provide an analogous classification in higher dimensions. For example, when $n\ge3$, the blow-up $X=\operatorname{Bl}_p\CP^n$ admits a K\"ahler metric with $H_g>0$ \cite{YZ}. Every smooth rational curve in $X$ has nonnegative total normal degree. But it is not a compact Hermitian symmetric space.

We also determine the optimal local pinching constants on $\FF_1$ and $\FF_2$. For $d\ge1$, let $\FF_d=\mathbb P\bigl(\OO_{\CP^1}\oplus\OO_{\CP^1}(-d)\bigr)$ and set
\[
\Lambda_{\FF_d}:=
\sup\{\lambda_{\FF_d}(g):g\text{ is a K\"ahler metric with }H_g>0\}.
\]
Alvarez, Chaturvedi, and Heier \cite{ACH} constructed Hitchin-type metrics with local pinching $(2d+1)^{-2}$ and showed this value is optimal within their family. Our bound for the negative section proves that their values for $d=1,2$ are also optimal without restricting the metric to that family.

\begin{corollary}\label{hirzebruch}
\[
\Lambda_{\FF_1}=\frac19,\qquad
\Lambda_{\FF_2}=\frac1{25},\qquad
\Lambda_{\FF_d}\le\frac1{25}\quad(d\ge3).
\]
\end{corollary}

\section{Pinching obstruction from negative normal degree}\label{obstruction}

The first lemma turns a negative orthogonal bisectional curvature component into a bound for the holomorphic sectional curvature ratio at a single point.

\begin{lemma}\label{two-plane}
Let $(M^n,g)$ be a K\"ahler manifold with positive holomorphic sectional curvature. If orthonormal vectors $u,w\in T_p^{1,0}M$ satisfy
\[
R(u,\bar u,w,\bar w)\le-\mu H_g(u)
\quad\text{for some }\mu>0,
\]
then $\lambda_{\{p\}}(g)\le\Phi(\mu)$.
\end{lemma}

\begin{proof}
After interchanging $u$ and $w$ if necessary, assume that
$\alpha:=H_g(u)\le\beta:=H_g(w)$. The assumed curvature inequality remains valid after this interchange. Put
$\gamma:=R(u,\bar u,w,\bar w)\le-\mu\alpha$.
For $0\le t\le1$, the vector
$\sqrt t\,u+e^{i\theta}\sqrt{1-t}\,w$ has unit length. Averaging over the phase gives
\begin{equation}\label{phase-average}
\begin{aligned}
H_{\min}(p)
&\le\frac1{2\pi}\int_0^{2\pi}
  H_g\bigl(\sqrt t\,u+e^{i\theta}\sqrt{1-t}\,w\bigr)\,d\theta\\
&=\alpha t^2+4\gamma t(1-t)+\beta(1-t)^2.
\end{aligned}
\end{equation}
Since $H_{\max}(p)\ge\beta$, set $\xi=\beta/\alpha\ge1$. Minimizing the resulting quadratic polynomial in $t$ at
$t_0=(\xi+2\mu)/(1+\xi+4\mu)\in(0,1)$ yields
\begin{equation}\label{pointwise-estimate}
\lambda_{\{p\}}(g)\le
\frac{\xi-4\mu^2}{\xi(1+\xi+4\mu)}=:f_\mu(\xi).
\end{equation}
A differentiation shows that $f_\mu$ increases up to $\xi_*=2\mu(1+4\mu)$ and decreases thereafter. Thus the maximum over $\xi\ge1$ occurs at $\xi=1$ if $\mu\le1/4$, and at $\xi=\xi_*$ if $\mu\ge1/4$. These two values are $1/2-\mu$ and $(1+4\mu)^{-2}$, respectively, as claimed.
\end{proof}

\begin{proof}[Proof of Theorem~\ref{main1}]
Set $r=n-1$. Let $e$ be a local unit tangent field along $C$, and choose local unitary normal fields $\nu_1,\ldots,\nu_r$. Write
\[
\alpha=H_g(e),\qquad
\Gamma=\sum_{j=1}^rR(e,\bar e,\nu_j,\bar\nu_j).
\]
Let $R^{TC}$ and $R^N$ denote the curvatures of the induced tangent and quotient Hermitian bundles, respectively. The Gauss and quotient curvature formulas give
\[
\alpha=R^{TC}(e,\bar e,e,\bar e)+|\Pi(e,e)|_g^2,
\qquad
\Gamma=\tr R^N(e,\bar e)-|\Pi(e,e)|_g^2.
\]
Since $\deg TC=2$ and $\deg N_{C/M}=-d$, Chern--Weil theory yields
\begin{equation}\label{integral-identities}
\int_C\alpha\,\omega_C=4\pi+S,
\qquad
\int_C\Gamma\,\omega_C=-2\pi d-S.
\end{equation}
For $\mu=\mu(n,d,S)$, these identities imply $\int_C(\Gamma+r\mu\alpha)\,\omega_C=0$. Hence there exists $p\in C$ with $\Gamma(p)\le-r\mu\alpha(p)$. The Hermitian form $R(e,\bar e,\cdot,\overline{\cdot})$ on the normal space at $p$ has trace $\Gamma(p)$, so some unit normal vector $\nu$ satisfies
\[
R(e,\bar e,\nu,\bar\nu)
\le\frac{\Gamma(p)}r\le-\mu H_g(e).
\]
Lemma~\ref{two-plane} now gives
$\lambda_C(g)\le\lambda_{\{p\}}(g)\le\Phi(\mu)$.
Finally,
\[
(n-1)\mu(n,d,S)
=\frac{4\pi}{4\pi+S}\,\frac d2
 +\frac{S}{4\pi+S}\cdot1
\ge\min\{d/2,1\}.
\]
Since $\Phi$ is decreasing, the remaining inequality in
\eqref{main-bound} follows.
\end{proof}

\begin{corollary}\label{minus-one}
Let $(M^2,g)$ be a K\"ahler surface with positive holomorphic sectional curvature and $\lambda_M(g)\ge1/9$. If $C\subset M$ is a smooth rational curve with $C^2=-1$, then $C$ is totally geodesic. Moreover, for every $p\in C$, every unit tangent vector $e\in T_p^{1,0}C$, and every unit normal vector $\nu\in(T_p^{1,0}C)^\perp$, one has
\[
R(e,\bar e,\nu,\bar\nu)=-\frac12H_g(e),\qquad
R(e,\bar e,e,\bar\nu)=0,\qquad
H_g(\nu)=H_{\max}(p)=3H_g(e).
\]
\end{corollary}

\begin{proof}
Along $C$, write
$\alpha=H_g(e)$, $\beta=H_g(\nu)$, and
$\gamma=R(e,\bar e,\nu,\bar\nu)$.
We first show that $\gamma+\alpha/2\ge0$ pointwise. Otherwise, at some point $\gamma<-\mu\alpha$ for a number $\mu>1/2$. Lemma~\ref{two-plane} would then give $\lambda_M(g)\le\Phi(\mu)<\Phi(1/2)=1/9$, a contradiction. Since $\deg N_{C/M}=C^2=-1$, the identities \eqref{integral-identities} imply
\[
0\le\int_C\left(\gamma+\frac12\alpha\right)\omega_C
=-\frac S2\le0.
\]
Consequently $S=0$, $C$ is totally geodesic and $\gamma=-\alpha/2$ along $C$.
Hence 
$R(e,\bar e,e,\bar\nu)=0$.

Fix $p\in C$. Taking $t=2/3$ in \eqref{phase-average} gives $H_{\min}(p)\le\beta/9$. But the pinching assumption also gives
\[
H_{\min}(p)\le\frac\beta9
\le\frac{H_{\max}(p)}9\le H_{\min}(p).
\]
Thus $\beta=H_{\max}(p)$ and $H_{\min}(p)=\beta/9$.
In particular, $\xi:=\beta/\alpha\ge1$. Applying
\eqref{pointwise-estimate} with $\mu=1/2$ yields
\[
\frac19=\lambda_{\{p\}}(g)
\le\frac{\xi-1}{\xi(\xi+3)}
=\frac19-\frac{(\xi-3)^2}{9\xi(\xi+3)}\le\frac19.
\]
Equality forces $\xi=3$, completing the proof.
\end{proof}

\section{Classification of K\"ahler surfaces}\label{classification}

\begin{proof}[Proof of Theorem~\ref{main2}]
Suppose $\lambda_M(g)\ge1/9$. By the contrapositive of
Lemma~\ref{two-plane}, for any orthonormal vectors
$e,\nu\in T_p^{1,0}M$ we have
\[
R(e,\bar e,\nu,\bar\nu)\ge-\frac12H_g(e).
\]
Since $\dim_{\mathbb C}M=2$, it follows that
\[
\operatorname{Ric}_g(e,\bar e)
=H_g(e)+R(e,\bar e,\nu,\bar\nu)
\ge\frac12H_g(e)>0.
\]
Thus $-K_M$ is positive and, by Kodaira's embedding theorem, ample. Consequently $M$ is a del Pezzo surface. The classification of del Pezzo surfaces \cite[Theorem~2.15]{MiroRoigPonsLlopis2014} says that $M$ is $\CP^2$, $\CP^1\times\CP^1$, or the blow-up of $\CP^2$ at $1\le k\le8$ points in general position. But every such blow-up contains a smooth rational $(-1)$-curve $C$ and then Theorem~\ref{main1} gives $\lambda_M(g)\le1/9$. Consequently, if \(\lambda_M(g)>1/9\), then \(M\) must be biholomorphic to either \(\mathbb{CP}^2\) or \(\mathbb{CP}^1\times\mathbb{CP}^1\). This proves (i).

We now show that if \(\lambda_M(g)=1/9\), then \(M\) cannot be the blow-up of \(\mathbb{CP}^2\) at more than one point.

Let $C$ be any
smooth rational $(-1)$-curve, fix $p\in C$, and choose orthonormal vectors
$e\in T_p^{1,0}C$ and $\nu\perp T_p^{1,0}C$. Set $\alpha=H_g(e)$.
By Corollary~\ref{minus-one},
\[
R(e,\bar e,\nu,\bar\nu)=-\frac\alpha2,
\qquad R(e,\bar e,e,\bar\nu)=0,
\qquad H_g(\nu)=3\alpha.
\]
The curvature inequality above shows that the function
\[
F(u):=\operatorname{Ric}_g(u,\bar u)-\frac12H_g(u)
\]
is nonnegative on the unit sphere in $T_p^{1,0}M$, and $F(e)=0$.
Consider the variation $u(s)=(e+s\nu)/\sqrt{1+s^2}$, $s \in \mathbb R$. Then
\[
\begin{aligned}
0
&=\left.\frac{d}{ds}\right|_{s=0}F(u(s))\\
&=2\operatorname{Re}\left(\operatorname{Ric}_g(\nu,\bar e)
  -R(\nu,\bar e,e,\bar e)\right)
=2\operatorname{Re}\operatorname{Ric}_g(\nu,\bar e).
\end{aligned}
\]
Replacing $\nu$ by $i\nu$ shows that $\operatorname{Ric}_g(e,\bar\nu)=0$. Therefore the Ricci tensor is diagonal in the unitary basis $\{e,\nu\}$, with eigenvalues
\[
\operatorname{Ric}_g(e,\bar e)=\frac\alpha2,
\qquad
\operatorname{Ric}_g(\nu,\bar\nu)=\frac{5\alpha}{2}.
\]
In particular, $\tr_g\operatorname{Ric}_g=3\alpha$.

Suppose another smooth rational $(-1)$-curve $C'$ passes through $p$.
Take a unit tangent vector $v\in T_p^{1,0}C'$ and a unit normal $v^\perp$. Applying Corollary~\ref{minus-one} to $C'$ yields
\[
\operatorname{Ric}_g(v,\bar v)=\frac12H_g(v),\qquad
\operatorname{Ric}_g(v^\perp,\overline{v^\perp})=\frac52H_g(v).
\]
Taking their sum and comparing with $\tr_g\operatorname{Ric}_g=3\alpha$, we obtain $H_g(v)=\alpha$. Thus $\operatorname{Ric}_g(v,\bar v)=\alpha/2$. Since this is the smallest Ricci eigenvalue at $p$, the vector $v$ lies in its one-dimensional eigenspace. Hence $T_p^{1,0}C'=T_p^{1,0}C$. We conclude that two distinct smooth rational $(-1)$-curves cannot meet transversely.

If $M$ were the blow-up of $\CP^2$ at $k\ge2$ points in general position, let $E_1$ be the exceptional curve over the first point and $L_{12}$ the strict transform of the line through the first two points. Then
\[
E_1^2=L_{12}^2=-1,\qquad E_1\cdot L_{12}=1.
\]
These are distinct $(-1)$-curves meeting transversely, a contradiction.
Therefore $k\le1$, and the possibilities in (ii) follow.
\end{proof}

\medskip
\noindent\textbf{Declaration on the use of AI.}
The mathematical results and proofs were developed by the author. ChatGPT was used primarily to help identify potential gaps in the proofs and to assist with language and exposition.

\noindent\textbf{Conflict of Interest.} The author declares no conflict of interest.

\noindent\textbf{Funding.} This research did not receive any specific grant from funding agencies in the public, commercial, or not-for-profit sectors.
\bibliographystyle{amsplain}
\bibliography{pinchingpaper}

\end{document}